\documentclass[12pt]{article} 
\usepackage{amsfonts,latexsym,color}
\input{psfig.sty}

\newcounter{environment}[section]
\renewcommand{\theenvironment}{%
\arabic{section}.\arabic{environment}}

{\begin{rm}\refstepcounter{environment}{\textbf\theenvironment\
\bf Definition.~~}}%
{\end{rm}}

{\begin{rm}\refstepcounter{environment}{\textbf\theenvironment\
\bf Example.~~}}%
{\end{rm}}

{\begin{rm}\refstepcounter{environment}{\textbf\theenvironment\
\bf Conjecture.~~}}%
{\end{rm}}

{\begin{rm}\refstepcounter{environment}{\textbf\theenvironment\
\bf Proposition.~~}}%
{\end{rm}}

{\begin{rm}\refstepcounter{environment}{\textbf\theenvironment\
\bf Proposition and Definition.~~}}%
{\end{rm}}

{\begin{rm}\refstepcounter{environment}{\textbf\theenvironment\
\bf Theorem.~~}}%
{\end{rm}}

{\begin{rm}\refstepcounter{environment}{\textbf\theenvironment\
\bf Corollary.~~}}%
{\end{rm}}

{\begin{rm}\refstepcounter{environment}{\textbf\theenvironment\
\bf Lemma.~~}}%
{\end{rm}}

\begin{document}
\newcommand{\qbc}[2]{{\left [{#1 \atop #2} \right ]}}
\newcommand{\df}[2]{1-x_{#1}x_{#2}^{-1}}
\newcommand{\pdf}[2]{\left(1-x_{#1}x_{#2}^{-1}\right)}
\newcommand{\beq}{\begin{equation}}
\newcommand{\eeq}{\end{equation}}
\newcommand{\be}{\begin{enumerate}}
\newcommand{\ee}{\end{enumerate}}
\newcommand{\bea}{\begin{eqnarray}}
\newcommand{\eea}{\end{eqnarray}}
\newcommand{\beas}{\begin{eqnarray*}}
\newcommand{\eeas}{\end{eqnarray*}}
\newcommand{\ds}{\displaystyle}
\newcommand{\bm}[1]{{\mbox{\boldmath $#1$}}}
\newcommand{\rr}{\mathbb{R}}
\newcommand{\zz}{\mathbb{Z}}
\newcommand{\nn}{\mathbb{N}}
\newcommand{\pp}{\mathbb{P}}
\newcommand{\st}{\,:\,}
\newcommand{\ps}{\mathbb{Q}[[x]]_0}
\newcommand{\cp}{{\cal P}}
\newcommand{\cf}{{\cal F}}
\newcommand{\gl}{\mathrm{GL}}
\newcommand{\sn}{\mathfrak{S}_n}
\newcommand{\cc}{\mathbb{C}}
\newcommand{\ch}{\mathrm{char}}
\newcommand{\vp}{\varphi}
\newcommand{\dt}{R^{(2)}}
\newcommand{\qtk}{K_{\lambda\mu}(q,t)}
\newcommand{\hs}{\mathrm{Hilb}^n(\cc^2)}
\newcommand{\is}{\mathrm{is}}
\newcommand{\tm}{\textcolor{magenta}}

\begin{centering}
\textcolor{red}{\Large\bf Recent Progress}\\
\textcolor{red}{\Large\bf in Algebraic Combinatorics}\\[.2in] 
\textcolor{blue}{Richard P. Stanley}\footnote{Partially supported by
  NSF grant 
\#DMS-9988459.}\\ 
Department of Mathematics\\
Massachusetts Institute of Technology\\
Cambridge, MA 02139\\
\emph{e-mail:} rstan@math.mit.edu\\[.1in]
\tm{\small version of 18 December 2001}\\[.2in]
\end{centering}
\section{Introduction.}
Algebraic combinatorics is alive and well at the dawn at the new
millenium. Algebraic combinatorics is difficult to define precisely;
roughly speaking it involves objects that can be interpreted both
combinatorially and algebraically, e.g., as the cardinality of a
combinatorially defined set and the dimension of an algebraically
defined vector space. Sometimes the combinatorial interpretation is
used to obtain an algebraic result, and sometimes \emph{vice versa}.
Mathematicians have been engaged in algebraic combinatorics at least
since Euler (in particular, his work on partitions), but it wasn't
until the 1960's, primarily under the influence of Gian-Carlo Rota,
that there was a systematic attempt to establish the foundations of
algebraic combinatorics and bring it into the mathematical mainstream.
This effort has been highly successful, and algebraic combinatorics
has by now become a mature and thriving discipline.

We have chosen three major breakthroughs to highlight recent work in
algebraic combinatorics. All three areas have initiated a flurry of
further work and suggest many further directions of research to keep
practitioners of algebraic combinatorics occupied well into the new
century. Our choice of topics was partially influenced by the relative
ease in describing the main results to nonexperts in algebraic
combinatorics.  Much other outstanding work has been done that is
not discussed here.

\section{The saturation conjecture.}
The saturation conjecture concerns certain integers known as
\emph{Littlewood-Richardson coefficients}. Given the theme of this
paper, it is not surprising that they have both an algebraic and a
combinatorial definition. First we discuss the algebraic definition,
which is more natural than the combinatorial one.

Let GL$(n,\cc)$ denote the group of all invertible transformations
from an $n$-dimensional complex vector space $V$ to itself. After
choosing an ordered basis for $V$ we may identify GL$(n,\cc)$ with the
group of $n\times n$ nonsingular matrices over the complex numbers
(with the operation of matrix multiplication). Consider the map
$\varphi:\gl(2,\cc)\rightarrow \gl(3,\cc)$ defined by
    $$ \varphi \left[\begin{array}{cc} a & b\\ c & d \end{array}\right] 
    =\left[ \begin{array}{ccc} a^2 & 2ab & b^2\\ 
     ac & ad+bc & bd\\ c^2 & 2cd & d^2
     \end{array} \right]. $$
This can be checked to be a group homomorphism (and hence a
\emph{representation} of $\gl(2,\cc)$ of degree 3). Moreover, the entries of
$\varphi(A)$ are polynomial functions of the entries of $A$. Hence
$\varphi$ is a \emph{polynomial representation} of $\gl(2,\cc)$. If
$A\in \gl(2,\cc)$ has eigenvalues $x,y$, then it can also be checked
that $\varphi(A)$ has eigenvalues $x^2,xy,y^2$. Define the
\emph{character} $\ch\,\varphi$ of $\varphi$ to be the trace
of $\varphi(A)$, regarded as a function of the eigenvalues $x,y$ of
$A$. Hence
   $$ \ch\,\varphi=x^2+xy+y^2. $$
It was first shown by Schur that the polynomial representations of
GL$(n,\cc)$ are \emph{completely reducible}, i.e., a direct sum of
irreducible representations. The nequivalent irreducible
polynomial representations $\varphi_\lambda$ of GL$(n,\cc)$ are
indexed by \emph{partitions} $\lambda=(\lambda_1,\dots,\lambda_n)$ of
length at most $n$, i.e., $\lambda_i\in\zz$ and
$\lambda_1\geq\cdots\geq \lambda_n\geq 0$. Moreover,
$\ch\,\varphi_\lambda$ is a symmetric function
$s_\lambda(x_1,\dots,x_n)$ that had been originally defined by Cauchy
and Jacobi and is now known as a \emph{Schur function}. A well-known
property of Schur functions is their \emph{stability}: 
  $$ s_\lambda(x_1,\dots,x_n,0) = s_\lambda(x_1,\dots,x_n). $$
For this reason we can let $n\rightarrow\infty$ and consider the Schur
function $s_\lambda$ in \emph{infinitely} many variables
$x_1,x_2,\dots$ and specialize to $x_1,\dots,x_n$ when dealing with
$\gl(n,\cc)$. For more
information on symmetric functions and the representation theory of
$\gl(n,\cc)$, see \cite{fulton}\cite{macd}\cite{ec2}.

If $A:V\rightarrow V$ and $B:W\rightarrow W$ are linear
transformations on finite-dimensional vector spaces, then
  $$ \mathrm{tr}(A\otimes B) = \mathrm{tr}(A)\cdot\mathrm{tr}(B), $$
where $A\otimes B$ denotes the tensor (or Kronecker) product of $A$
and $B$, acting on $V\otimes W$. 
Hence if $\lambda$, $\mu$, and $\nu$ are partitions and we set
  $$ c_{\mu\nu}^\lambda=\mathrm{mult}(\varphi_\lambda,\varphi_\mu
         \otimes \vp_\nu), $$
the multiplicity of $\vp_\lambda$ in the tensor product $\vp_\mu
\otimes \vp_\nu$ (when written as a direct sum of irreducible
representations), then
  $$ s_\mu s_\nu = \sum_\lambda c_{\mu\nu}^\lambda s_\lambda. $$

The nonnegative integers $c_{\mu\nu}^\lambda$ are known as
\emph{Littlewood-Richardson coefficients}, and the
\emph{Littlewood-Richardson rule}
\cite[Ch.\ 5]{fulton}\cite[{\S}I.9]{macd}\cite[Appendix~A1.3]{ec2} gives a
combinatorial interpretation of them (which we will not state
here). If $m$ is a positive integer and
$\lambda=(\lambda_1,\lambda_2,\dots)$ a partition, then write 
$m\lambda = (m\lambda_1,m\lambda_2,\dots)$.

\textbf{Saturation conjecture.} If $c_{m\mu,m\nu}^{m\lambda}\neq 0$,
then $c_{\mu\nu}^\lambda\neq 0$.

The saturation conjecture was proved recently by Allen Knutson and
Terence Tao \cite{k-t}\cite{k-t2} using a new \emph{honeycomb model}
for describing Littlewood-Richardson coefficients. An elegant
exposition of the proof was given by Anders Buch \cite{buch}, and a
detailed survey of all the material in this section (and more) was
given by William Fulton \cite{fulton2}. A proof of the saturation
conjecture based on representations of quivers was later given by Harm
Derksen and Jerzy Weyman \cite{d-w}.

Why is the proof of the saturation conjecture an important
breakthrough? The answer is that it is related in a surprising way to
a number of other topics. The first concerns the
eigenvalues of hermitian matrices. Let $A,B,C$ be
$n\times n$ hermitian matrices. Hence their eigenvalues are
real. Denote the eigenvalues of $A$ as
    $$ \alpha\,:\ \alpha_1\geq\cdots\geq \alpha_n, $$
and similarly $\beta$ and $\gamma$ for $B$ and $C$. 
Considerable attention has been given to the following problem.

\textbf{Problem.} Characterize those triples $(\alpha,\beta,\gamma)$ for
which there exist hermitian matrices $A+B=C$ with eigenvalues
$\alpha$, $\beta$, and $\gamma$.

By taking traces we see that
  \beq  \sum \gamma_i = \sum \alpha_i+\sum \beta_i. \label{eq:tr} \eeq
After much work by a number of researchers, A. Horn conjectured a
complete characterization of triples $(\alpha,\beta,\gamma)$,
consisting of (\ref{eq:tr}) together with linear inequalities of the
form
  \beq \sum_{k\in K}\gamma_k\leq \ds\sum_{i\in I}\alpha_i +
       \ds\sum_{j\in J}\beta_j, \label{eq:horn} \eeq
for certain sets
   $$ I,J,K\subset \{1,\dots,n\},\ \ \ |I|=|J|=|K|. $$
For instance, when $n=2$ Horn's inequalities (which are easy to show
that together with (\ref{eq:tr})
characterize $(\alpha,\beta,\gamma)$ in this case) become
  \beas \gamma_1 & \leq & \alpha_1+\beta_1\\
         \gamma_2 & \leq & \alpha_2+\beta_1\\
         \gamma_2 & \leq & \alpha_1+\beta_2. \eeas
For $n=3$ there are twelve inequalities, as follows:
  \beas \gamma_1 & \leq & \alpha_1+\beta_1\\
       \gamma_2 & \leq & \min(\alpha_1+\beta_2,\alpha_2+\beta_1)\\
      \gamma_3 & \leq & \min(\alpha_1+\beta_3,\alpha_2+\beta_2,
            \alpha_3+\beta_1)\\
      \gamma_1+\gamma_2 & \leq & \alpha_1+\alpha_2+\beta_1+\beta_2\\
      \gamma_1+\gamma_3 & \leq & \min(\alpha_1+\alpha_2+\beta_1
          +\beta_3,\alpha_1+\alpha_3+\beta_1+\beta_2)\\
     \gamma_2+\gamma_3 & \leq & \min(\alpha_1+\alpha_2+\beta_2+\beta_3,
         \alpha_1+\alpha_3+\beta_1+\beta_3,\alpha_2+
         \alpha_3+\beta_1+\beta_2).  \eeas

The connection between the Saturation Conjecture and Horn's conjecture
was given by Alexander Klyachko \cite{kly}.

\textbf{Theorem.} \emph{The Saturation Conjecture implies Horn's
  conjecture.} 

A more precise connection between Littlewood-Richardson coefficients
and eigenvalues of hermitian matrices is provided by the following
result, implicit in the work of Heckman \cite{heck} and more
explicit in Klyachko \cite{kly}. 

\textbf{Theorem.} \emph{Let $\alpha,\beta$, and $\gamma$ be partitions
  of length at most $n$. The Saturation Conjecture implies that the
  following two conditions are equivalent:}
\vspace{-.2in}
\begin{itemize} \item $c_{\alpha\beta}^\gamma\neq 0$.
   \item  \emph{There exist $n\times n$ hermitian matrices $A+B=C$
       with eigenvalues $\alpha,\beta$, and $\gamma$.}
  \end{itemize}
Since equation (\ref{eq:horn}) consists of \emph{linear} inequalities,
the two theorems above show that the nonvanishing of
$c_{\alpha\beta}^\gamma$ depends on (explicit) linear inequalities
among the coordinates of $\alpha, \beta, \gamma$. Thus for fixed $n$
the points $(\alpha,\beta,\gamma)\in\rr^{3n}$ for which
$c_{\alpha\beta}^\gamma\neq 0$ are the integer points in a certain
convex cone. Hence the subject of polyhedral combinatorics is
closely associated with the theory of Littlewood-Richardson
coefficients. For further information on this point of view, see
\cite{zel}. 

The theorems stated above involve hermitian matrices. It is known
\cite[Thm.\ 3]{fulton2} that exactly the same results hold for the
class of real symmetric matrices.

There are a number of other situations in which Littlewood-Richardson
coefficients play a surprising role. These situations are thoroughly
discussed in \cite{fulton2}. We mention one of them here. Given a
partition $\lambda =(\lambda_1,\lambda_2,\dots)$ and a prime $p$, let
$G$ be a (finite) abelian $p$-group of type $\lambda$, i.e.,
  $$ G\cong \left(\zz/p^{\lambda_1}\zz\right)\times 
        \left(\zz/p^{\lambda_2}\zz\right)\times\cdots. $$
Given further partitions $\mu$ and $\nu$, let $g_{\mu\nu}^\lambda(p)$
denote the number of subgroups $H$ of $G$ of type $\mu$ such that the
quotient group $G/H$ has type $\nu$. 

\textbf{Theorem.} (a) \emph{$g_{\mu\nu}^\lambda(p)$ is a polynomial
  function of $p$  with integer coefficients.} 
 
 (b) \emph{For any prime $p$ we have that $g_{\mu\nu}^\lambda(p)\neq
  0$ if and only if $c_{\mu\nu}^\lambda\neq 0$.}

The polynomial $g_{\mu\nu}^\lambda(t)$ is called a \emph{Hall
  polynomial} after the pioneering work of Philip Hall \cite{hall}.
Hall established the above theorem, except that in part (b) he only
showed that $g_{\mu\nu}^\lambda(t)$ vanishes identically (as a
polynomial in $t$) if and only if $c_{\mu\nu}^\lambda= 0$.
Subsequently Miller Maley \cite{maley} showed that the polynomial
$g_{\mu\nu}^\lambda(t+1)$ has nonnegative coefficients, from which (b)
follows. For an exposition of the basic properties of Hall
polynomials, see \cite[Chs.\ II and III.2]{macd}. The theory of Hall
polynomials holds in the more general context of the ring of integers
(i.e., the unique maximal order) of a division algebra of finite rank
over a $\mathfrak{p}$-adic field \cite[Remark 3, p.\ 179]{macd} or
even more generally for $q$-primary lattices \cite[Thm.\ 
4.81]{tesler}.

\section{The $n!$ and $(n+1)^{n-1}$ conjectures.} \label{sec:n}
The $n!$ and $(n+1)^{n-1}$ conjectures concern the action of the
symmetric group $\sn$ on two sets $(x_1,\dots,x_n)$ and
$(y_1,\dots,y_n)$ of $n$ variables. In order to appreciate these
conjectures, knowledge of the situation for one set of $n$ variables
is of value. We therefore first review this theory (for which the
proofs are \emph{much} easier). $\sn$ acts on the polynomial ring
$A=\cc[x_1,\dots,x_n]$ by permuting variables, i.e., for $w\in\sn$ let
$w\cdot x_i=x_{w(i)}$ and extend to all of $A$ in the obvious way.
Let
  $$ A^{\sn}=\{ f\in A\st w\cdot f=f\ \ \ \forall w\in \sn\}, $$
the \emph{ring of invariants} of the action of $\sn$ on $A$. The
invariant polynomials $f\in A^{\sn}$ are the \emph{symmetric
  polynomials} in the variables $x_1,\dots,x_n$ (over $\cc$). The
``fundamental theorem of symmetric functions'' asserts that 
  $$ A^{\sn}=\cc[e_1,\dots,e_n], $$
a polynomial ring in the algebraically independent elementary
symmetric functions
  $$ e_k = \sum_{1\leq i_1<\cdots<i_k\leq n}
         x_{i_1}\cdots x_{i_k}. $$
Regard $n$ as fixed and define the ring
  $$ R = A/(e_1,\dots,e_n). $$
The ring $R$ inherits the usual grading from $A$, i.e.,
  $$ R=R_0\oplus R_1\oplus\cdots, $$
where $R_i$ is spanned by (the images of) all homogeneous
polynomials of degree $i$ in the variables $x_1,\dots,x_n$. Because
the generators $e_1,\dots,e_n$ of $R^{\sn}$ are algebraically
independent of degrees $1,2,\dots,n$, it is easy to see that 
  $$ \dim_\cc R=n!, $$
and more generally,
  \beq \sum_i \dim_\cc\!\left( R_i\right)q^i =
  (1+q)(1+q+q^2)\cdots(1+q+\cdots+q^{n-1}), \label{eq:qnfact} \eeq
the standard ``$q$-analogue'' of $n!$.

Since the ideal $(e_1,\dots,e_n)$ of $R$ is $\sn$-invariant, $\sn$
acts on $R$. Moreover, this action respects the grading of
$R$, i.e., $w\cdot R_i=R_i$ for all $w\in\sn$. Thus
$R$ is in fact a \emph{graded $\sn$-module}, and we can ask, as a
refinement of (\ref{eq:qnfact}), for the multiplicity of each
irreducible representation of $\sn$ in $R_i$. For the action on
$R$ as a whole the situation is simple to describe (and not
difficult to prove): $R$ affords the \emph{regular
representation} of $\sn$, i.e., the multiplicity of each irreducible
representation is its degree (or dimension).

To describe the $\sn$-module structure of $R_i$, we need some
understanding of the (inequivalent)
irreducible representations of $\sn$. They are
indexed by partitions $\lambda$ of $n$ (denoted $\lambda\vdash n$),
i.e, $\lambda= (\lambda_1,\dots,\lambda_\ell)\in\zz^\ell$ where
$\lambda_1\geq \cdots\geq \lambda_\ell>0$ and $\sum \lambda_i=n$. The
dimension of the irreducible $\sn$-module $M_\lambda$ indexed by
$\lambda\vdash n$ is denoted by
$f^\lambda$ and is equal to the number of \emph{standard Young
  tableaux} (SYT) of shape $\lambda$, i.e., the number of ways to
insert the numbers $1,2,\dots,n$ (without repetition) into an array
of shape $\lambda$ (i.e., left-justified with $\lambda_i$ entries in
row $i$) so that every row and column is increasing. For instance
$f^{(3,2)}=5$, as shown by the five SYT
  $$ \begin{array}{lclclclcl}
     1\,2\,3 & & 1\,2\,4 & & 1\,2\,5 & & 1\,3\,4 & & 1\,3\,5\\
     4\,5  & & 3\,5  & & 3\,4  & & 2\,5  & & 2\,4 \end{array}. $$
There is also a simple explicit formula (e.g.,
\cite[Exam.\ I.5.2]{macd}\cite[Cor.\ 7.21.6]{ec2}), known as the
\emph{hook-length formula}, for $f^\lambda$.

Since $R$ affords the regular representation of $\sn$, the
multiplicity of $M_\lambda$ in $R$ is equal to $f^\lambda$. Thus
we would like to describe the multiplicity of $M_\lambda$ in
$R_i$ as the number of SYT $T$ of shape $\lambda$ with some
additional property depending on $i$. This property is the value of
the \emph{major index} of $T$, denoted MAJ$(T)$. It is defined by
  $$ \mathrm{MAJ}(T)=\ds\sum_{i+1\, \mathrm{below}\, i\, 
    \mathrm{in}\, T} i, $$  
where the sum ranges over all entries $i$ of $T$ such that $i+1$
appears in a lower row than $i$. For instance, the SYT of shape
$(3,2,2)$ shown below has MAJ$(T)=2+3+6=11$.
  $$ T=\begin{array}{l} 1\,\bm{2}\,\bm{6}\\
           \bm{3}\,5\\ 4\,7 \end{array}. $$
The following result is due independently to Lusztig (unpublished) and
Stanley \cite[Prop.\ 4.11]{rs:inv}. 

\textbf{Theorem.} \emph{Let $\lambda\vdash n$. Then}
  $$ \mathrm{mult}(M_\lambda,R_i) = \#\{ \mathrm{SYT}\ T\st
    \mathrm{shape}(T)=\lambda,\ \mathrm{MAJ}(T)=i \}. $$

For example, let $n=5$. There are three SYT with five entries and
major index 3, namely,
  $$ \begin{array}{lcclccl} 1\,2\,\bm{3}\,5 & & & 1\,2\,\bm{3} & & &
           \bm{1}\,4\,5\\ 4 & & & 4\,5 & & & \bm{2}\\
            & & & & & & 3 \end{array}. $$
It follows that 
  $$ R_3\cong M_{41}\oplus M_{32}\oplus M_{311}. $$
There is another description of $R$ which leads to a
different generalization to two sets of $n$ variables. Given any
polynomial $P(x_1,\dots,x_n)$ over $\cc$, define $\partial P$ to be
the complex vector space spanned by $P$ and all its partial
derivatives of all orders. For instance $\partial(x+y)^2$ has
dimension four, one basis being $\{(x+y)^2,x,y,1\}$. Let
  \beq V_n=\prod_{1\leq i<j\leq n}(x_i-x_j). \label{eq:vdm} \eeq
It is easy to see that
  $$ R \cong \partial V_n $$
as graded $\sn$-modules. In particular, $\dim(\partial V_n)=n!$ and
$\partial V_n$ affords the regular representation of $\sn$.

Adriano Garsia and Mark Haiman had the idea of generalizing the above
constructions of $R$ and $\partial V_n$ to two sets
$x=(x_1,\dots,x_n)$ and $y=(y_1,\dots,y_n)$ of $n$ variables. 
For the first generalization, let $\sn$ act \emph{diagonally} on
$B=\cc[x,y]$, i.e.,
  $$ w\cdot x_i =x_{w(i)},\ \ w\cdot y_i=y_{w(i)}. $$
Let 
  $$ B^{\sn} =\{f\in B\st w\cdot f=f\ \ \ \forall w\in\sn\}, $$
the ring of invariants of the action of $\sn$ on $B$. It is no
longer the case that $B^{\sn}$ is generated by algebraic independent
elements. (For general information about rings of invariants of finite
groups, see for instance \cite{smith}\cite{rs:inv}.) However, we can
still define 
  $$ \dt = S/I, $$ 
where $I$ is the ideal of $B$ generated by elements of $B^{\sn}$
with zero constant term. The $(n+1)^{n-1}$ conjecture of Garsia and
Haiman \cite{g-h}\cite{g-h2} was recently proved
by Haiman \cite{hai:npo}, based on techniques he
developed to prove the $n!$ conjecture discussed below, together
with a theorem of Bridgeland, King, and Reid on the McKay
correspondence.

\textbf{Theorem} ($(n+1)^{n-1}$ conjecture). $\dim_\cc
\dt=(n+1)^{n-1}$ 

Just as $R$ had the additional structure of a graded
$\sn$-module, similarly $\dt$ is a \emph{bigraded}
$\sn$-module. In other words, 
  $$ \dt = \bigoplus_{i,j}\dt_{ij}\ \ \mbox{(vector
  space direct sum)}, $$
where $R_{ij}^{(2)}$ is the subspace of $\dt$ spanned by
(the images of) polynomials that are homogeneous of degree $i$ in the
$x$ variables and degree $j$ in the $y$ variables, and moreover
$\dt_{ij}$ is invariant under the action of $\sn$ on
$\dt$. For instance, when $n=4$ it can be computed that
  $$ \dt_{2,1}\cong 2M_{211}\oplus M_{22}\oplus M_{31}. $$
In particular, 
 $$ \dim_\cc  \dt_{2,1}=2f^{211}+f^{22}+f^{31} =
      2\cdot 3+2+3=12. $$
Garsia and Haiman stated in \cite{g-h3} (see also \cite[Conj.\
7.5]{haiman1}) a complicated conjectured formula for
mult$(M_\lambda,\dt_{ij})$. Haiman's proof of the
$(n+1)^{n-1}$ conjecture mentioned above actually establishes this
stronger conjecture of Garsia and Haiman.
A consequence of Haiman's result asserts the following
\cite{g-h3}\cite[p.\ 246]{haiman1}. Let $\Gamma$ be the
anti-invariant subspace of $\dt$, i.e.,
  $$ \Gamma = \{f\in\dt\st w\cdot f = \mathrm{sgn}(w)f\ \ \ \forall
       f\in\sn\}, $$
where sgn$(w)$ denotes the sign of the permutation $w$. Then
   $$ \dim_\cc  \Gamma = \ds\frac{1}{n+1}{2n\choose n}, $$
a \emph{Catalan number}. James Haglund \cite{hag} conjectured and
Garsia and Haglund \cite{g-hag} proved a combinatorial interpretation of
the $\Gamma$ bigrading, i.e., a combinatorial interpretation of the
numbers $\dim_\cc \Gamma_{ij}$. For some information on the ubiquitious
appearance of Catalan (and related) numbers throughout mathematics,
see \cite[Exer.\ 6.19--6.38]{ec2} and the addendum at
www-math.mit.edu/$\sim$rstan/ec.html.

The number $\dim_\cc\dt=(n+1)^{n-1}$ has a number of combinatorial
interpretations, e.g., it is the number of forests of rooted trees on
$n$ vertices \cite[Prop.\ 5.3.2]{ec2} or the number of parking
functions of length $n$ \cite[Exer.\ 5.49]{ec2}. It is natural to ask
whether one can give a combinatorial interpretation of $\dim_\cc  \dt_{ij}$
that refines some known interpretation of $(n+1)^{n-1}$. At present
this question is open.

We turn to the second generalization of $R$ due to to Garsia and
Haiman. First we need to define a generalization of the Vandermonde
product (\ref{eq:vdm}) to two sets of variables. Let $\mu\vdash
n$. Coordinatize the squares of the diagram of $\mu$ by letting
$(i-1,j-1)$ be the coordinate of the square in the $i$th row and $j$th
column. For instance, the coordinates of the squares of the diagram of
$\mu=(3,2)$ are given by

\vspace{.3in}
\centerline{\psfig{figure=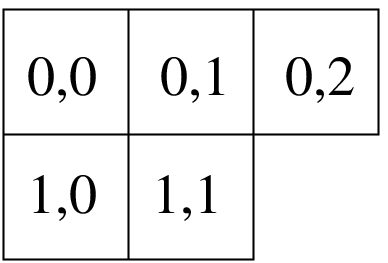}}
  
Let $(i_1,j_1),\dots,(i_n,j_n)$ be the coordinates of the squares of
the diagram of $\mu$ (in some order), and define the $n\times n$
determinant 
  $$ D_\mu = \left| x_r^{i_s}y_r^{j_s}\right|_{r,s=1,\dots,n}. $$
For instance,
  $$ D_{32}=\left| \begin{array} {lllll} 1 & y_1 & y_1^2 & x_1 & x_1y_1 \\
      1 & y_2 & y_2^2 & x_2 & x_2y_2 \\
      1 & y_3 & y_3^2 & x_3 & x_3y_3 \\
      1 & y_4 & y_4^2 & x_4 & x_4y_4 \\
      1 & y_5 & y_5^2 & x_5 & x_5y_5 \end{array} \right|. $$
Note that if $\mu$ consists of a single row (i.e., $\mu$
consists of the single part $n$) then $D_\mu=V_n(y)$, while if
$\mu$ consists of a single column then $D_\mu=V_n(x)$. 

The $n!$ conjecture of Garsia and Haiman \cite{g-h}\cite{g-h2}, later
proved by Haiman \cite{haiman:proof}, is the following
assertion.

\textbf{Theorem} ($n!$ conjecture). \emph{For any $\mu\vdash n$, we
  have} $$ \dim_\cc  \partial D_\mu = n!. $$

  The space $\partial D_\mu$, just as $\dt$, is a bigraded
  $\sn$-module. For each $i,j\geq 0$ and $\lambda\vdash n$, we can ask
  for a ``description'' of the integer mult$\left(M_\lambda,
    \left(D_\mu\right)_{ij}\right)$. Garsia and Haiman
  \cite{g-h}\cite{g-h2} gave such a description, and Haiman
  \cite[Thm.\ 5.4]{haiman1} showed that it actually followed from the 
  $n!$ conjecture. The Garsia-Haiman description involves the theory of
  \emph{Macdonald symmetric functions}, a generalization of Schur
  functions due to I. G. Macdonald \cite{macdsf}\cite[Ch.\ VI]{macd}
  and currently of great interest in several different areas, such as
  the representation theory of quantum groups, affine Hecke algebras,
  and the Calegero-Sutherland model in particle physics (see
  \cite{haiman:proof} for references). We won't define Macdonald
  symmetric functions here but will give a brief indication of
  Haiman's result.

Let $\lambda,\mu\vdash n$. The coefficient of $x^\mu=x_1^{\mu_1}
x_2^{\mu_2}\cdots$ in the Schur function $s_\lambda$ is known as a
\emph{Kostka number}, denoted $K_{\lambda\mu}$, and has a simple
combinatorial interpretation in terms of semistandard Young tableaux
\cite[(5.13)]{macd}\cite[{\S}7.10]{ec2}. In the theory of Macdonald
polynomials there arises naturally a two-parameter generalization
$\qtk$ of the Kostka number $K_{\lambda\mu}=
K_{\lambda\mu}(0,1)$. \emph{A priori} $\qtk$ is only a
rational function of $q$ and $t$, but Macdonald conjectured that it
was a polynomial with nonnegative integer coefficients. In 1996--98
several independent proofs were given that $\qtk$ was indeed a
polynomial with integer coefficients, but nonnegativity remained open.
Haiman showed the remarkable fact that $\qtk$ is essentially the
bigraded Hilbert series for the $\lambda$-isotypic component of
$D_\mu$. More precisely,
  $$ t^{b(\mu)} K_{\lambda\mu}(q,1/t) = \sum_{r,s\geq 0}
    \mathrm{mult}\left(M_\lambda,
   \left( D_\mu\right)_{r,s}\right)\,t^rq^s, $$
where $b(\mu)=\sum (i-1)\mu_i$. This formula establishes the
nonnegativity of the coefficients of $\qtk$, though a combinatorial
interpretation of these coefficients remains open.

Hamian's proof is based on the geometry of the Hilbert scheme
$\hs$ of $n$ points in the plane. (Claudio Procesi suggested
to Haiman the possible relevance of the Hilbert scheme.) Let $X$ and
$Y$ be indeterminates. We can define
$\hs$ as a set by
  $$ \hs = \{ I\subseteq \cc[X,Y]\st
        \dim_\cc\cc[X,Y]/I=n\}, $$
i.e., all ideals $I$ of $\cc[X,Y]$ such that the quotient ring
$\cc[X,Y]/I$ is an $n$-dimensional vector space. Suppose that ${\cal
Z}=\{z_1,\dots,z_n\}$ is a set of $n$ \emph{distinct} points in
$\cc^2$. Let
  $$ I_{{\cal Z}}=\{ f\in\cc[X,Y]\st f(z_1)=\cdots=f(z_n)=0\}. $$
Then $I_{{\cal Z}}$ is an ideal of $\cc[X,Y]$ such that
$\cc[X,Y]/I_{{\cal Z}}$ can be identified with the space of all
functions $f:{\cal Z}\rightarrow \cc$, so $I_{{\cal Z}}\in
\hs$. This explains why $\hs$ is called the Hilbert scheme of $n$
points in the plane --- it is a closure of the space of all
$n$-element subsets of $\cc^2$. In fact, $\hs$ has the structure of a
smooth irreducible algebraic variety, of dimension
$2n$. 

The remarkable connections between $\hs$ and the $n!$ and
$(n+1)^{n-1}$ conjectures are too technical to discuss here, but let
us give a vague hint or two. Write $H^n=\hs$.
Given a partition $\mu\vdash n$, let 
$U_\mu$ be the set of all ideals $I\in H^n$ such that a basis for
$\cc[x,y]/I$ consists of the (images of the) monomials $x^hy^k$, where
the $(h,k)$'s are the coordinates for the squares of the diagram of
$\mu$. Then the sets $U_\mu$ are open, affine, and cover $H^n$,
suggesting the possible relevance of $H^n$ to the $n!$
conjecture. Moreover, for each $I\in H^n$ there is a natural way to
associate an $n$-element \emph{multiset} $\pi(I)\subset \cc^2$. The
$n$-element multisets contained in $\cc^2$ form an affine variety
Sym$^n(\cc^2)$, viz.,
  $$ \mathrm{Sym}^n(\cc^2) = (\cc^2)^n/\sn = \mathrm{Spec}\ 
   \cc[x_1,\dots,x_n,y_1,\dots,y_n]^{\sn}, $$
suggesting the possible relevance of $H^n$ to the
$(n+1)^{n-1}$ conjecture. See the papers \cite{haiman1} and
\cite{haiman:proof} for details. 

It is natural to ask about generalizing the work of Garsia and Haiman
to more than two sets of variables. However, all obvious conjectures
turn out to be false. One difficult is that the Hilbert scheme
Hilb$^n(\cc^k)$ is no longer smooth for $k>2$.

The $(n+1)^{n-1}$ and $n!$ conjectures are just the beginning of an
amazing edifice of conjectures due to Garsia, Haiman, and their
collaborators. For instance, we defined a determinant $D_\lambda$ when
$\lambda$ is a partition of $n$, regarded as a certain subset of
$\nn\times\nn$ (where $\nn=\{0,1,2,\dots\}$). In exactly the same way
we can define $D_X$ for \emph{any} $n$-element subset $X$ of
$\nn\times\nn$. Bergeron, Garsia, and Tesler \cite{b-g-t} then
conjecture (and prove in some special cases) for several classes of
subsets $X$ that $\dim_\cc (\partial D_X)=k_Xn!$ for some positive
integer $k_X$; and in fact $\partial D_X$, regarded as an
$\sn$-module, affords $k_X$ copies of the regular representation.

\section{Longest increasing subsequences.}
Let $w=a_1a_2\cdots a_n\in\sn$. An \emph{increasing subsequence} of $w$
is a subsequence $a_{i_1}a_{i_2}\cdots a_{i_k}$ of $w$ for which
$a_{i_1} <a_{i_2}<\cdots<a_{i_k}$. Let $\is_n(w)$ denote the length of
the longest increasing subsequence of $w\in\sn$. For instance, if
$w=274163958\in\mathfrak{S}_9$ then $\is_9(w)=4$, exemplified by the
increasing subsequences 2469 and 1358. There has been much recent
interest in the behavior of the function $\is_n(w)$. A survey of much
of this work has been given by Percy Deift \cite{deift}. 

The first question of interest is the expected value $E(n)$ of
$\is_n(w)$, where $w$ ranges uniformly over $\sn$. Thus
  $$ E(n) = \frac{1}{n!}\sum_{w\in\sn}\is_n(w). $$
Elementary arguments show that 
  $$ \frac 12\sqrt{n}\leq E(n) \leq e\sqrt{n}, $$
and Hammersley \cite[Thm.\ 4]{ham} showed in 1972, using subadditive
ergodic theory, that the limit
  $$ c =\lim_{n\rightarrow\infty} \frac{E(n)}{\sqrt{n}} $$
exists. Vershik and Kerov \cite{v-k} (with the difficult direction
$c\geq 2$ shown independently by Logan and Shepp \cite{l-s}) showed in
1977 that $c=2$.

The proof of Vershik-Kerov and Logan-Shepp is based on the identity
  \beq E(n) = \frac{1}{n!}\sum_{\lambda\vdash n}\lambda_1
  \left( f^\lambda\right)^2, \label{eq:exen} \eeq
where $\lambda=(\lambda_1,\lambda_2,\dots)$ and $f^\lambda$ denotes
the number of SYT of shape $\lambda$ as in
Section~\ref{sec:n}. Equation (\ref{eq:exen}) is due to Craige
Schensted \cite{schen} and is an immediate consequence of the
Robinson-Schensted-Knuth algorithm; see also \cite[Exer.\
7.109(a)]{ec2}. 

The work of Vershik-Kerov and Logan-Shepp only determines the
asymptotic behavior of the expectation of is$_n(w)$. What about
stronger results? A major breakthrough was made by Jinho Baik, Percy
Deift, and Kurt Johansson \cite{b-d-j}, and has inspired much
further work. To describe their results, let Ai$(x)$ denote the
\emph{Airy function}, viz., the unique solution to the second-order
differential equation 
  $$ \mathrm{Ai}''(x) = x\,\mathrm{Ai}(x), $$
subject to the condition
  $$ \mathrm{Ai}(x)\sim \frac{e^{-\frac 23 x^{3/2}}}
       {2\sqrt{\pi}x^{1/4}}\ \mathrm{as}\ x\rightarrow \infty. $$
Let $u(x)$ denote the unique solution to the nonlinear third order
equation
  \beq u''(x) = 2u(x)^3+xu(x), \label{eq:pain} \eeq
subject to the condition
 $$ u(x) \sim -\mathrm{Ai}(x),\ \mathrm{as}\ x\rightarrow\infty. $$
Equation (\ref{eq:pain}) is known as the \emph{Painlev\'e II
  equation}, after Paul Painlev\'e (1863--1933)\footnote{In addition
  to being a distinguished mathematician, in 1908 Painlev\'e was the
  first passenger of Wilbur Wright, during which they set a flight
  duration record of 70 minutes, and in 1917 and 1925 he held a
  position equivalent to Prime Minister of France.}. Painlev\'e
completely classified differential equations (from a certain class of
second order equations) whose ``bad'' singularities (branch points and
essential singularities) were independent of the initial conditions.
Most of the equations in this class were already known, but a few were
new, including equation (\ref{eq:pain}).

Now define the \emph{Tracy-Widom distribution} to be the probability
distribution on $\rr$ given by 
  \beq F(t)= \exp\left( -\int_t^\infty (x-t)u(x)^2\,dx\right).
   \label{eq:t-w} \eeq
It is easily seen that $F(t)$ is indeed a probability distribution,
i.e., $F(t)\geq 0$ and $\int_{-\infty}^\infty F(t)dt=1$. Let $\chi$ be
a random variable with distribution $F$, and let $\chi_n$ be the
random variable on $\sn$ defined by
  $$ \chi_n(w) =\frac{\is_n(w)-2\sqrt{n}}{n^{1/6}}. $$
We can now state the remarkable results of Baik, Deift, and
Johansson. 

\textbf{Theorem.} \emph{As $n\rightarrow\infty$, we have}
  $$ \chi_n\rightarrow \chi\quad\mbox{in distribution}, $$
\emph{i.e., for all $t\in \rr$,}
  $$ \lim_{n\rightarrow\infty} \mathrm{Prob}(\chi_n\leq t) =F(t). $$

\textbf{Theorem.} \emph{For any $m=0,1,2,\dots$,}
  $$ \lim_{n\rightarrow\infty} E(\chi_n^m)=E(\chi^m). $$

\textbf{Corollary.} \emph{We have}
  \beas \lim_{n\rightarrow\infty}\frac{\mathrm{Var}(\is_n)}
     {n^{1/3}} & = &  \int t^2\,dF(t) -\left(\int
       t\,dF(t)\right)^2\\ & = & 0.8132\cdots, \eeas
\emph{where} Var \emph{denotes variance, and}
   \begin{eqnarray} \lim_{n\rightarrow\infty}\frac{E(\is_n)-2\sqrt{n}}
      {n^{1/6}} & = & \int t\,dF(t) \label{eq:e2} \\
     & = & -1.7711\cdots. \nonumber \end{eqnarray}
The above theorems are a vast refinement of the
  Vershik-Kerov and Logan-Shepp results concerning $E(n)$, the
  expectation of is$_n(w)$. The first theorem gives the entire
  limiting distribution (as $n\rightarrow\infty$) of is$_n(w)$, while
  the second theorem gives an asymptotic formula for the $m$th
  moment. Note that equation (\ref{eq:e2}) may be rewritten
  $$ E(n) = 2\sqrt{n}+\alpha n^{1/6}+o\left(n^{1/6}\right), $$
where $\alpha=\int t\,dF(t)$, thereby giving the second term in the
asymptotic behavior of $E(n)$.

We will say only a brief word on the proof of the above results,
explaining how combinatorics enters into the picture. Some kind of
analytic expression is needed for the distribution of is$_n(w)$. Such
an expression is provided by the following result of Ira Gessel
\cite{gessel}, later proved in other ways by various persons.

\textbf{Theorem.} \emph{Let} 
  \beas u_k(n) & = & \#\{ w\in\sn\st\is_n(w)\leq k\}\\
    U_k(x) & = & \sum_{n\geq 0} u_k(n)\frac{x^{2n}}{n!^2}\\
    B_i(x) & = & \sum_{n\geq 0}\frac{x^{2n+i}}{n!\,(n+i)!}. \eeas
\emph{Then}
  $$ U_k(x) = \det\left( B_{|i-j|}(x)\right)_{i,j=1}^k. $$

\textbf{Example.} We have
  \beas U_2(x) & = & \left| \begin{array}{cc} B_0(x) & B_1(x)\\
                B_1(x) & B_0(x) \end{array} \right|\\
           & = &  B_0(x)^2-B_1(x)^2. \eeas
From this it is easy to deduce that 
  $$ u_2(n) = \frac{1}{n+1}{2n\choose n}, $$
a Catalan number. This result was first stated by John Michael
Hammersley in 1972, with the first published proofs by Knuth
\cite[{\S}5.1.4]{knuth} and Rotem \cite{rotem}. There is a
more complicated expression for $u_3(n)$ due to Gessel
\cite[{\S}7]{gessel}\cite[Exer.\ 7.16(e)]{ec2}, namely,
  $$ u_3(n) =\frac{1}{(n+1)^2(n+2)} \sum_{j=0}^n {2j\choose j}
    {n+1\choose j+1}{n+2\choose j+2}, $$
while no ``nice'' formula for $u_k(n)$ is known for fixed $k>3$. 

Gessel's theorem reduces the theorems of Baik, Deift, and Johansson to
``just'' analysis, viz., the Riemann-Hilbert problem in the theory of
integrable systems, followed by the method of steepest descent to
analyze the asymptotic behavior of integrable systems. For further
information see the survey \cite{deift} of Deift mentioned above.

The asymptotic behavior of is$_n(w)$ (suitably scaled) turned out to be
identical to the Tracy-Widom distribution $F(t)$ of equation
(\ref{eq:t-w}). It is natural to ask how the Tracy-Widom distribution
arose in the first place. It seems surprising that such an
``unnatural'' looking function as $F(t)$ could have arisen
independently in two different contexts. Originally the Tracy-Widom
distribution arose in connection with the \emph{Gaussian Unitary
  Ensemble} (GUE). GUE is a certain natural probability distribution
on the space of all $n\times n$ hermitian matrices $M=(M_{ij})$,
namely,
  $$ Z_n^{-1} e^{-\mathrm{tr}(M^2)}dM, $$
where $Z_n$ is a normalization constant and
  $$ dM = \prod_i dM_{ii}\cdot\prod_{i<j} d(\mathrm{Re}\, M_{ij})
   d(\mathrm{Im}\,M_{ij}). $$
Let the eigenvalues of $M$ be $\alpha_1\geq \alpha_2\geq \cdots \geq
\alpha_n$. The following result marked the eponymous appearance
\cite{t-w} of the Tracy-Widom distribution:
  \beq \lim_{n\rightarrow\infty}\mathrm{Prob}\left( \left(\alpha_1-
    \sqrt{2n}\right)\sqrt{2}n^{1/6}\leq t\right) = F(t). \label{eq:tw}
   \eeq
Thus as $n\rightarrow\infty$, is$_n(w)$ and $\alpha_1$ have the same
distribution (after scaling). 

It is natural to ask, firstly, whether there is a result analogous to
equation (\ref{eq:tw}) for the other eigenvalues $\alpha_k$ of the GUE
matrix $M$, and, secondly, whether there is some connection between
such a result and the behavior of 
increasing subsequences of random permutations. A generalization of
(\ref{eq:tw}) was given by Tracy and Widom \cite{t-w} (expressed in
terms of the Painlev\'e II function $u(x)$). The connection with
increasing subsequences was conjectured in \cite{b-d-j} and proved
independently by Borodin-Okounkov-Olshanski \cite{b-o-o}, Johannson
\cite{johan}, and Okounkov \cite{ok}. Given $w\in\sn$, define integers
$\lambda_1,\lambda_2,\dots$ by letting $\lambda_1+\cdots+\lambda_k$ be
the largest number of elements in the union of $k$ increasing
subsequences of $w$. For instance, let $w=247951368$. The longest
increasing subsequence is 24568, so $\lambda_1=5$. The largest union
of two increasing subsequences is 24791368 (the union of 2479 and
1368), so $\lambda_1+\lambda_2=8$. (Note that it is impossible to find
a union of length 8 of two increasing subsequences that contains an
increasing subsequence of length $\lambda_1=5$.) Finally $w$ itself is
the union of the three increasing subsequences 2479, 1368, and 5, so
$\lambda_1+\lambda_2+\lambda_3=9$. Hence
$(\lambda_1,\lambda_2,\lambda_3)=(5,3,1)$ (and $\lambda_i=0$ for
$i>3$). Readers familiar with the theory of the
Robinson-Schensted-Knuth algorithm will recognize the sequence
$(\lambda_1, \lambda_2, \dots)$ as the \emph{shape} of the two
standard Young tableaux obtained by applying this algorithm to $w$, a
well-known result of Curtis Greene \cite{greene}\cite[Thm.\ 
A1.1.1]{ec2}. (In particular, $\lambda_1\geq \lambda_2\geq\cdots$, a
fact which is by no means obvious.)  The result of
\cite{b-o-o}\cite{johan}\cite{ok} asserts that as as
$n\rightarrow\infty$, $\lambda_k$ and $\alpha_k$ are equidistributed,
up to scaling.

The Tracy-Widom distribution arose completely independently in the
behaviour of is$_n(w)$ and GUE matrices. Is this connection just a
coincidence? The work of Okounkov \cite{ok} provides a connection,
\emph{via} the theory of random topologies on surfaces.

\end{document}